\documentclass[leqno]{amsart}
\usepackage{amsmath}
\usepackage{amssymb}
\usepackage{amsthm}
\usepackage{enumerate} 
\usepackage[mathscr]{eucal}

\newtheorem{theorem}{Theorem}[section]

\newtheorem{example}[theorem]{Example}

\newtheorem{remark}[theorem]{Remark}

\usepackage{cite} 
\usepackage{url}  
\usepackage{ifthen}  
\usepackage{multicol}   
\usepackage[utf8]{inputenc} 
\newboolean{publ}

\begin{document}

		\title[On strong orthogonality and  strictly convex normed linear space]{On strong orthogonality and  strictly convex normed linear spaces}
	\author[K.Paul, D. Sain and K. Jha]{Kallol Paul , Debmalya Sain and Kanhaiya Jha}
	
		\address[Paul]{Department of Mathematics, Jadavpur University, Kolkata 700032, India.}
	\email{kalloldada@gmail.com}
	
	\address[Sain]{Department of Mathematics, Jadavpur University, Kolkata 700032, India.}
	\email{saindebmalya@gmail.com}

	\address[Jha]{Department of Mathematical Sciences, School of Science, Kathmandu University, POBox Number 6250, Kathmandu, NEPAL.}
	\email{jhakn@ku.edu.np}

	\begin{abstract}
		We introduce the notion of strongly orthogonal set relative to an element  in the sense of Birkhoff-James in a normed linear space to 
		find a necessary and sufficient condition for an element $ x $ of the unit sphere $ S_{X}$ to be an exposed point of the unit ball $ B_X .$ We then 
		prove that a normed linear space is strictly convex iff  for each element x of the unit sphere there exists 
		a bounded linear operator A on X which attains its norm only at the points of the form  $ \lambda x $ with $ \lambda \in S_{K} $. 
		
	\end{abstract}
	
	\subjclass[2020]{Primary 46B20, Secondary 47A30}
	\keywords{Orthogonality, Strict convexity, Extreme point.}

	\maketitle

		\section{Introduction.}
		\noindent Suppose $ (X, \|.\|) $ is a normed linear space over the field K, real or complex. X is said to be strictly convex iff every element of the unit sphere $ S_X = \{ x \in X : \|x \| = 1 \} $ is an extreme point of the unit ball $ B_X = \{ x \in X : \|x \| \leq 1 \} .$ There are many equivalent characterizations of the strict convexity of a normed space, some of them given in [8,12] are\\
		(i) If $ x,y \in S_X $ then we have $ \| x + y \| < 2,$\\
		(ii) Every non-zero continuous linear functional attains a maximum on atmost one point of the unit sphere, \\
		(ii) If $ \| x+y\| = \|x\| + \|y\|,~ x \neq 0 $ then $ y = c x $ for some $ c \geq 0. $ \\
		The notion of  strict convexity plays an important role in the studies of the geometry of Banach Spaces. One may go through  [1,3-5,6-8,11-14] for more information related to strictly convex spaces. \\
		\noindent An element $x$ is said to be orthogonal to $y$ in X in the sense of Birkhoff-James [2,7,8], written as, $x \bot_{B}y $,  iff 
		\[ \| x \| \leq \| x + \lambda y \| ~~ \mbox{for all  scalars} ~\lambda .\]
		\noindent If X is an inner product space then $ x \bot_{B}y $ implies $ \| x \| < \| x + \lambda y \| $ for all  scalars $ \lambda \neq 0.$ Motivated by this fact we here introduce the notion of strong orthogonality as follows:\\
		\noindent \textbf{Strongly orthogonal in the sense of Birkhoff-James}: In a normed linear space X an element $x$ is said to be strongly orthogonal to another element $y$  in the sense of Birkhoff-James, written as $ x \perp_{SB}y,  $  iff 
		\[ \| x \| < \| x + \lambda y \| ~~ \mbox{for all  scalars} ~\lambda \neq 0 .\]
		If $ x \perp_{SB}y  $ then $ x \perp_{B}y  $ but the converse is not true. In $ l_\infty ({R}^2) $ the element $ (1,0) $ is orthogonal to $(0,1)$ in the sense of Birkhoff-James but not strongly orthogonal.\\
		\textbf{Strongly orthogonal  set relative to an element} : A finite set of elements $ S = \{x_1,x_2, \ldots x_n \} $ is said to be a strongly orthogonal set  relative to an element $ x_{i_0}$ contained in S in the sense of Birkhoff-James    iff 
		\[  \| x_{i_0} \| < \| x_{i_0} + \sum_{j=1,j\neq i_0}^{n} \lambda_j x_j \|\]
		whenever not all $ \lambda_j$'s are 0. \\
		An infinite set of elements is said to be a strongly orthogonal set relative to an element contained in the set in the sense of Birkhoff-James  iff every finite subset containing that element is strongly orthogonal relative to that element in the sense of Birkhoff-James.\\
		\textbf{Strongly orthogonal Set}:  A finite set of elements $ \{x_1,x_2, \ldots x_n \} $ is said to be a strongly orthogonal set in the sense of Birkhoff-James   iff for each  $ i \in  \{ 1,2, \ldots n \} $
		\[  \| x_{i} \| < \| x_i + \sum_{j=1,j\neq i}^{n} \lambda_j x_j \|\]
		whenever not all $ \lambda_j$'s  are 0. \\
		An infinite set of elements is said to be a  strongly orthogonal set in the sense of Birkhoff-James   iff every finite subset of the set is a strongly  orthogonal set in the sense of Birkhoff-James. \\
		\noindent Clearly if  a set is strongly orthogonal in the sense of Birkhoff-James then it is strongly orthogonal relative to every element of the set in the sense of Birkhoff-James. If X has a Hamel basis which is strongly orthogonal in the sense of Birkhoff-James then we call the Hamel basis a strongly  orthogonal Hamel  basis in the sense of Birkhoff-James and if X has a Hamel basis which is strongly orthogonal relative to an element of the basis  in the sense of Birkhoff-James then we call the Hamel basis a strongly  orthogonal Hamel basis relative to that element of the basis  in the sense of Birkhoff-James. If in addition the norm of each element of a strongly orthogonal set is 1 then accordingly we call them orthonormal. \\
		As for example  the set $ \{(1,0,\ldots ,0), (0,1,0 \ldots ,0), \ldots ,(0,0,\ldots,1) \} $ is a strongly orthonormal Hamel basis in the sense of Birkhoff-James in $ l_{1}({R}^n) $ but not in $ l_{\infty}({R}^n).$ \\
		In $ l_2({R}^3)$ the set $ \{ (1,0,0),(0,1,0),(0,1,1) \} $ is strongly orthogonal relative to (1,0,0) in the sense of Birkhoff-James but not relative to(0,1,1). \\
		\noindent In this paper we give another characterization of strictly convex normed linear spaces by using the Hahn-Banach theorem and  the notion of strongly orthogonal Hamel  basis relative to an element in the sense of Birkhoff-James, more precisely we explore the relation between the existence of strongly orthogonal Hamel basis relative to an element with unit norm  in the sense of Birkhoff-James in a normed space and that of an extreme point of the unit ball in the space. We also prove that a  normed linear space is strictly convex iff for each point x of the unit sphere there exists a bounded linear operator A on X which attains its norm only at the points of the form  $ \lambda x $ with $ \lambda \in S_K.$
		
		\section{Main Results}
		
		\noindent We first obtain a sufficient condition for an element in the unit sphere to be an extreme point of the unit ball in  an arbitrary normed linear space.
		\begin{theorem}
			Let X be a normed linear space and $ x_0 \in S_X $. If there  exists a Hamel basis of X containing $ x_0 $ which is strongly orthonormal  relative to $ x_0$ in the sense of Birkhoff-James then $ x_0 $ is an extreme point of $ B_X.$
		\end{theorem}
		\noindent \textbf{Proof.} Let D = $ \{ x_0, x_{\alpha} : \alpha \in \Lambda \} $ be a strongly  orthonormal Hamel basis  relative to $ x_0$  in the sense of Birkhoff-James. \\
		If possible suppose $ x_0 $ is not an extreme point of $ B_X $, then  $ x_0 = t z_1 + (1-t)z_2 $ where $ 0<t<1 $ and $ \|z_1\|= \|z_2\|=1. $\\ 
		So there exists $ \alpha_1, \alpha_2, \ldots \alpha_n $ in $ \Lambda $ such that 
		\[  z_1 = \beta_0 x_0 + \sum_{j=1}^{n} \beta_j x_{\alpha_j} ~\mbox{and}~ z_2 = \gamma_0 x_0 + \sum_{j=1}^{n} \gamma_j x_{\alpha_j} \]
		for some scalars $\beta_j, \gamma_j(j=0,1,2,\dots n).$  \\
		If $\beta_0 = 0 $ and $ \gamma_0 = 0 $ then $ x_0 = t z_1 + (1-t)z_2 $ implies that 
		\[ x_0 = \sum_{j=1}^{n}(t\beta_{j} + (1-t)\gamma_j)x_{\alpha_{j}}, \] 
		which contradicts the fact that every finite subset of $D$ is linearly independent. So the case $\beta_0=0 $ and $ \gamma_0=0 $ is ruled out. \\
		If $\beta_0 \neq 0,~~ \gamma_0 = 0 $ then  as $ \{ x_0, x_{\alpha_1}, x_{\alpha_2}, \ldots , x_{\alpha_n} \} $  is a strongly  orthonormal set  relative to $ x_0$  in the sense of Birkhoff-James so we get 
		\[ 1 = \|z_1\| = | \beta_0| \| x_0 + \sum_{j=1}^{n} \frac{\beta_j}{\beta_0} x_{\alpha_{j}} \| \geq | \beta_0 |. \]
		Now 
		\[ x_0 = t \beta_0 x_0 + \sum_{j=1}^{n} (t \beta_j + (1-t)\gamma_j) x_{\alpha_j} \]
		and so $ t\beta_0 = 1 $ which is not possible as $|\beta_0| \leq 1 $ and $ 0 < t < 1$.\\
		\noindent Similarly $\beta_0 = 0, ~~ \gamma_0 \neq 0$ is also ruled out. \\
		Thus we have $\beta_0 \neq 0 $ and $  \gamma_0 \neq 0.$ \\
		Our claim is that at least one of  $ \mid \beta_0 \mid , \mid \gamma_0 \mid $ must be less than 1. \\
		If possible suppose $ \mid \beta_0 \mid > 1 .$ Then 
		\[ \| \beta_0 x_0 + \sum_{j=1}^{n} \beta_j x_{\alpha_j} \| = \mid \beta_0 \mid \| x_0 + \sum_{j=1}^{n} \frac{\beta_j}{\beta_0} x_{\alpha_j} \| \geq \mid \beta_0 \mid > 1 . \]
		This contradicts that $ \| z_1 \| = 1 .$  Thus $ \mid \beta_0 \mid \leq 1 .$ Similarly  $ \mid \gamma_0 \mid \leq 1 .$ 
		We next show that $ \mid \beta_0 \mid = 1 $ and  $ \mid \gamma_0 \mid = 1 $ can not hold simultaneously.\\
		\textbf{Case1}. X is a real normed linear space. \\
		Then $ \mid \beta_0 \mid = 1 $ implies that 
		\[ 1 = \| z_1\| = \mid \beta_0 \mid \| x_0 + \sum_{j=1}^{n} \frac{\beta_j}{\beta_0} x_{\alpha_j} \| > \| x_0 \|, \]
		unless $ \beta_i = 0~~ \forall i=1,2, \ldots n. $ \\
		Thus $ \mid \beta_0 \mid = 1 \Rightarrow z_1 = \beta_0 x_0 \Rightarrow z_1= \pm x_0 \Rightarrow x_0=z_1=z_2~ or~ t=0,$ which is not possible. Thus $ \mid \beta_0 \mid \neq 1.$  Similarly $ \mid \gamma_0 \mid \neq 1 $.\\
		\textbf{Case2}. X is a complex normed linear space. \\
		Then $ \mid \beta_0 \mid = 1 $ implies that 
		\[ 1 = \| z_1\| = \mid \beta_0 \mid \| x_0 + \sum_{j=1}^{n} \frac{\beta_j}{\beta_0} x_{\alpha_j} \| > \| x_0 \|, \]
		unless $ \beta_i = 0~~ \forall i=1,2, \ldots n. $ \\
		Thus $ \mid \beta_0 \mid = 1 \Rightarrow z_1 = \beta_0 x_0 \Rightarrow z_1= e^{i \theta} x_0 $, similarly $ \mid \gamma_0 \mid = 1 \Rightarrow z_2 = e^{i \phi} x_0 $. Then $ x_0 = t e^{i \theta}x_0 + (1-t) e^{i \phi} x_0 \Rightarrow x_0=z_1=z_2,$ which is not possible.
		Thus $ \mid \beta_0 \mid = 1 $ and  $ \mid \gamma_0 \mid = 1 $ can not hold simultaneously.\\
		So at least one of  $ \mid \beta_0 \mid , \mid \gamma_0 \mid $ is less than 1.\\
		Now $ x_0 = t z_1 + (1-t)z_2 $ implies 
		\[ t \beta_0 + (1-t)\gamma_0 = 1, ~  t \beta_j + (1-t)\gamma_j = 0~ \forall ~j = 1,2, \ldots n. \]
		But $ \mid \beta_0 \mid < 1 $ or $  \mid \gamma_0 \mid < 1$  implies 
		\[ 1=  \mid t \beta_0 + (1-t)\gamma_0 \mid  \leq t \mid \beta_0 \mid + (1-t) \mid \gamma_0 \mid < 1, \]
		which is not possible.\\
		Thus $ x_0 $ is an extreme point of $ B_X.$ This completes the proof.\\
		\noindent The converse of the above Theorem  is however not always true. If $x_0$ is an extreme point of $ B_X $ then there may or  may not exist a strongly orthonormal Hamel basis relative to $x_0$ in the sense of Birkhoff-James.
		\begin{example}
			(i) Consider $ ({R}^2,\|.\|) $ where the unit sphere S is given by
			$ S = \{ (x,y) \in {R}^2 : x = \pm 1 ~ and ~ -1 \leq y \leq 1 \} \cup \{ (x,y) \in {R}^2 : x^2-2y+y^2=0 ~ and~ y \geq 1 \} \cup \{ (x,y) \in {R}^2 : x^2+2y+y^2=0 ~ and ~ y \leq - 1 \}.$ Then $ (1,1) $ is an extreme point of the unit ball but there exists no strongly orthonormal Hamel basis relative to $(1,1)$ in the sense of Birkhoff-James. \\
			\noindent (ii) Consider $ ({R}^2,\|.\|) $ where the unit sphere S is given by
			$ S = \{ (x,y) \in {R}^2 : x = \pm 1 ~ and ~ -1 \leq y \leq 1 \} \cup \{ (x,y) \in {R}^2 : x^2+2y-3=0 ~ and~ y \geq 1 \} \cup \{ (x,y) \in {R}^2 : x^2-2y-3=0 ~ and ~ y \leq - 1 \}.$ Then $ (1,1) $ is an extreme point of the unit ball and $ \{ (1,1),(-1,1) \} $ is a strongly orthonormal basis relative to (1,1) in the sense of Birkhoff-James.\\
			\noindent (iii) In $ l_\infty({R}^3)$ the  extreme points of the unit ball are of the form $ (\pm1,\pm1,\pm1)$ and   for the extreme  point (1,1,1) we can find a strongly  orthonormal basis relative to (1,1,1) in the sense of Birkhoff-James which is $ \{ (1,1,1),(1,0,-1),(0,1,-1) \}. $\\
		\end{example}
		\noindent In the first two examples the extreme point (1,1) is such that every neighbourhood of (1,1) contains both extreme as well as non-extreme points whereas in the third case the extreme point (1,1,1) is an isolated extreme point.\\
		\noindent An element x in the boundary of a convex set S is called an \textbf{exposed point} of S iff there exists a hyperplane of support H to S through x such that $ H \cap S = \{x\}. $ The notion of exposed points can be found in [4,9,10,15]. We next prove that if the extreme point $ x_0 $ is an exposed point of $ B_X$ then there  exists a Hamel basis of X containing $ x_0 $ which is strongly orthonormal  relative to $ x_0$ in the sense of Birkhoff-James. 
		\begin{theorem}
			Let X be a normed linear space and $ x_0 \in S_X $ be an exposed point of $ B_X$. Then there  exists a Hamel basis of X containing $ x_0 $ which is strongly orthonormal  relative to $ x_0$ in the sense of Birkhoff-James.
		\end{theorem}
		\noindent \textbf{Proof.} As $ x_0$ is an exposed point of $ B_X$ so there exists a hyperplane of support H to $ B_X $ through $x_0$ such that $ H \cap B_X = \{x_0\} .$ Then we can find a linear functional f on X such that $ H = \{ x \in X : f(x) = 1 \}.$ Let $ H_0 = \{ x \in X : f(x) = 0 \}.$ Then $ H_0 $ is a subspace of X. Let $ D = \{ x_\alpha : \alpha \in \Lambda \} $ be a Hamel basis of $ H_0 $ with $ \| x_\alpha \| = 1 .$ Clearly $ \{ x_0 \} \cup D $ is a Hamel basis of X. We claim that $ \{ x_0 \} \cup D $ is a strongly  orthonormal set relative to $  x_0 $ in the sense of Birkhoff-James. \\
		Consider a finite subset $ \{ x_{\alpha_1}, x_{\alpha_2}, \ldots ,x_{\alpha_{n-1}} \} $ of D  and let $ ( \lambda_1, \lambda_2, \ldots ,\lambda_{n-1}) \neq ( 0,0, \ldots ,0 ).$ Now if $ z = x_0 + \sum_{j=1}^{n-1} \lambda_j x_{\alpha_j} $ then
		\begin{eqnarray*}
			& & f(z) = f(x_0 + \sum_{j=1}^{n-1} \lambda_j x_{\alpha_j}) = f(x_0) = 1 \\
			& & \Rightarrow z \in H \\
			& & \Rightarrow z \notin B_X, ~~as ~~ H \cap B_X = \{x_0\} 
		\end{eqnarray*}
		So $ \| x_0 + \sum_{j=1}^{n-1} \lambda_j x_{\alpha_j} \| > 1 = \| x_0 \|. $
		Thus $ \{ x_0 \} \cup D $ is a Hamel  basis containing $ x_0 $ which is strongly orthonormal   relative to $ x_0 $ in the sense of Birkhoff-James.\\ 
		This completes the proof.\\
		\noindent We next prove that 
		\begin{theorem}
			Let X be a  normed linear space and $ x_0 \in S_X $.  If there exists a Hamel basis of X containing $ x_0 $ which is strongly orthonormal  relative to $ x_0$ in the sense of Birkhoff-James then there exists  a bounded invertible linear operator A on X such that $ \| A \| = \| A x_0 \| > \| A y \| $ for all y in $ S_X $ with   $ y \neq \lambda x_0 $, $ \lambda \in S_K $.
		\end{theorem}
		\noindent \textbf{Proof.} Let $ \{ x_0, x_\alpha : \alpha \in \Lambda \} $ be a Hamel basis of X which is strongly orthonormal relative to $ x_0 $ in the sense of Birkhoff-James.\\
		Define a linear operator A on X by $ A(x_0) = x_0 $ and $ A(x_\alpha) = \frac{1}{2} x_{\alpha}~ \forall \alpha \in \Lambda .$ \\
		Clearly  A is invertible.
		Take any $ z \in X $ such that $ \| z \| = 1 .$  Then  $ z = \lambda_0 x_0 + \sum_{j=1}^{n-1} \lambda_j x_{\alpha_j} $ for some scalars $ \lambda_j$'s  and $ \lambda_0$.\\
		\noindent If  $ \lambda_0 = 0 $ then  $ Az = \frac{1}{2} z $ and so 
		\[ \| Ax_0 \| = 1 > \frac{1}{2} = \|Az\|.\]
		\noindent If $ \lambda_0 \neq 0 $ then as $  \{x_{0},x_{\alpha}:\alpha \in \Lambda\}$ is a strongly orthonormal Hamel basis relative to $ x_{0} $ in the sense of Birkhoff-James so we get 
		\[ 1 = \|z\| = \| \lambda_0 x_0 + \sum_{j=1}^{n-1} \lambda_j x_{\alpha_j} \| \geq | \lambda_0|.\]
		Hence we get 
		\begin{eqnarray*}
			\|Az \| & = & \| \lambda_0 x_0 + \frac{1}{2} \sum_{j=1}^{n-1} \lambda_j x_{\alpha_j} \| \\
			& = & \| \frac{1}{2}(\lambda_0 x_0 + \sum_{j=1}^{n-1} \lambda_j x_{\alpha_j}) + \frac{1}{2} \lambda_0 x_0 \| \\
			& \leq & \frac{1}{2} \|z\| + \frac{1}{2} | \lambda_0|\\
			& \leq & 1 = \|Ax_0\|.
		\end{eqnarray*}
		This proves that $ \|A\| \leq 1. $ Also $ \|Az\| = 1 $ iff $ \mid \lambda_0 \mid = 1 $ and $ \lambda_{j} = 0 ~ \forall j = 1,2, \ldots n-1.$ \\
		Thus $ \| Az \| = 1 $ iff $  z = \lambda_0 x_0 $ with $ \lambda_0 \in S_K.$  
		This completes the proof.\\
		\noindent We now prove that 
		\begin{theorem}
			Let X be a normed linear space and $ x_0 \in S_X$. If there exists a bounded linear operator $ A : X \rightarrow X $ which attains its norm only at the points of the form $ \lambda x_0 $ with $ \lambda \in S_K $ then $ x_0 $ is an exposed point of $ B_X.$
		\end{theorem}
		\noindent \textbf{Proof.} Assume without loss of generality that $ \| A \| = 1 $ and by the Hahn-Banach theorem there exists $ f \in S_{X^{*}} $  such that $ f(A(x_0)) = 1. $   \noindent Clearly $ \| foA \| = 1 $ as $ f \in S_{X^*} $ and $ \| A \| = \|Ax_0\| = 1. $ If $ y \in S_X $ is such that $ \mid foA(y) \mid = 1,$ then $ \| Ay\| = 1. $ \\ 
		Now  $ \| A\| = 1 $ and $A$ attains its norm only at the points of the form $ \lambda x_0 $ with $ \lambda \in S_K,$  so $ y \in \{ \lambda x_0 : \lambda \in S_K \}. $ \\
		Thus $ f o A $ attains its norm only at the points of the form $ \lambda x_0 $ with $ \lambda \in S_K .$  Considering the hyperplane 
		$ H = \{ x \in X : foA(x) = 1 \} $ it is easy to verify that $ H \cap B_X = \{x_0\} $ and so $ x_0 $ is an exposed point of $ B_X.$ \\
		
		\noindent Thus we obtained complete characterizations of exposed points which is stated clearly in the following theorem
		\begin{theorem}
			For a normed linear space X and a point $ x \in S_X$ the following are equivalent : \\
			1. x is an exposed point of $ B_X$.\\
			2. There exists a Hamel basis of X containing $ x $ which is strongly orthonormal  relative to $ x$ in the sense of Birkhoff-James.\\
			3. There exists a bounded linear operator A on X which attains its norm only at the points of the form  $ \lambda x $  with $ \lambda \in S_K $.
		\end{theorem}
		\noindent We  next give a characterization of strictly convex space as follows : 
		\begin{theorem}
			For a normed linear space X  the following are equivalent : \\
			1. X is strictly convex.\\
			2. For each $ x \in S_X $ there exists a Hamel basis of X containing $ x $ which is strongly orthonormal  relative to $ x$ in the sense of Birkhoff-James.\\
			3. For each $ x \in S_X $ there exists a bounded linear operator A on X which attains its norm only at the points of the form  $ \lambda x $  with $ \lambda \in S_K $.
		\end{theorem}
		\noindent \textbf{Proof.} The proof follows from previous theorem and the fact that a normed linear space $ X$  is strictly convex iff every  element of $ S_X$ is an exposed point of $B_X$.
		\begin{remark}
			Even though the notions of strong Birkhoff-James orthogonality and Birkhoff-James orthogonality coincide in Hilbert space they do not characterize Hilbert spaces as $ (R^n, \|.\|_p)( 1 < p < \infty, p \neq 2)$ is not a Hilbert space but the notions of strong Birkhoff-James orthogonality and Birkhoff-James orthogonality coincide there.
		\end{remark}
		
		\section*{Competing interests}
		The authors declare that they have no competing interests.

		\section*{Acknowledgements}
		We would  like to thank the referee for their invaluable suggestion. We would also like to thank Professor T. K. Mukherjee for his invaluable suggestion while preparing this paper.
		The first author would like to thank Jadavpur University and DST, Govt. of India  for the partial financial support provided through DST-PURSE project and the second author would like to thank UGC, Govt. of India for the financial support. \\

		
		
	
\end{document}